\documentclass[12pt,leqno]{amsart}
\usepackage{amssymb,latexsym}
\usepackage{enumerate}
\usepackage{amsfonts}
\usepackage{amsmath}
\makeatletter
\@namedef{subjclassname@2010}{%
  \textup{2010} Mathematics Subject Classification}
\newtheorem{thm}{Theorem}[section]

\newtheorem{lem}[thm]{Lemma}

\begin{document}

\title{Distribution of coefficients of modular forms and the
partition function}

\author[Shi-chao Chen]{Shi-chao Chen}
\address{Institute of Contemporary Mathematics, Department of Mathematics and Information Sciences, Henan
University, Kaifeng, 475004, P. R. China} \email{schen@henu.edu.cn}

\begin{abstract} Let $\ell\ge5$ be
an odd prime and $j, s$ be positive integers. We study the
distribution of the coefficients of integer and half-integral weight
modular forms modulo odd positive integer $M$. As a consequence,  we
prove that for each integer $1\le r\le\ell^j$,
$$\sharp\{1\le n\le X\ |\ p(n)\equiv r\pmod{\ell^j}\}\gg_{s,r,\ell^j}\frac{\sqrt X}{\log X}(\log\log X)^s.$$
\end{abstract}



\maketitle

\section{Introduction and results}

Let $p(n)$ be the partition function, that is, the number of ways to
write the positive integer $n$ as the sum of a non-increasing
sequence of positive integers. M. Newman \cite{ne} conjectured that
if $M$ is a positive integer, then for every integer $r$ there are
infinitely many non-nonegative $n$ such that $p(n)\equiv r(\text{
mod }M).$ This conjecture has been settled for many cases. By the
work of Atkin, Kolberg, Newman, and Kl\o ve \cite{at,kol,ne,kl}, it
is known that the conjecture is ture for $M=2,5,7,13,17,19,29,31$.
Ono and Ahlgren \cite{ao,a} obtained an algorithm which presumably
proves the truth of the conjecture for any given $M$ coprime to
$2\cdot3\cdot5\cdot7\cdot11$. Combining the results of Bruinier and
Ono [6, 7], Ahlgren and Boylan [2] show that the conjecture is true
for all primes $\ell\ge5$. Recently, they can prove that the
conjecture holds for all prime powers $\ell^j$ [3]. More precisely,
they obtained the following bounds:
\begin{align}\sharp\{1\le n\le X
| p(n)\equiv r(\text{ mod }\ell^j)\}\gg_{r,\ell^j}
\left\{\begin{array}{ll}\frac{\sqrt{X}}{\log X} \quad &\text{ if
}r\not\equiv0\pmod{\ell^j},\\ X\quad&\text{ if
}r\equiv0\pmod{\ell^j}.\end{array}\right.\end{align}

In this paper we will use Serre's observations on Galois
representations to obtain some distribution properties of the
coefficients of integral and half-integral weight modular forms. As
a consequence, we can slightly improve the bound (1) when
$r\not\equiv0(\text{ mod }\ell^j)$.

Before we state main results, we recall some facts on modular forms
(see \cite{ko} for reference). For integers $k\ge0, N\ge1$, let
$S_k(\Gamma_0(N),\chi)$ be the space of cusp forms of weight $k$
with respect to $\Gamma_0(N)$ with character $\chi$, where $\chi$ is
a Dirichlet character modulo $N$. Suppose that $\lambda\ge0$ is an
integer and $4|N$. We denote by
$S_{\lambda+\frac{1}{2}}(\Gamma_0(N),\chi)$ the space of cusps of
weight $\lambda+\frac{1}{2}$ with respect to congruence subgroup
$\Gamma_0(N)$ with character $\chi$.

Let $M$ be a positive integer and $f(z)=\sum_{n=1}^\infty a(n)q^n\in
S_{\lambda+\frac{1}{2}}(\Gamma_0(N),\chi)\cap \mathbb{Z}[[q]]$. We
call that the coefficients of $f(z)$ are {\it well-distributed}
modulo $M$ if, for any integer $r$ and positive integer $s$, we have
\begin{align*}\sharp\{1\le n\le X | a(n)&\equiv r\pmod M\}\\
&\gg_{f,s,r,M} \left\{\begin{array}{ll}\frac{\sqrt{X}}{\log
X}(\log\log X)^s \quad &\text{ if }r\not\equiv0\pmod M,\\
X\quad&\text{ if }r\equiv0 \pmod M.\end{array}\right.\end{align*}

\noindent{\it Remark 1.1.} Our definition of the well-distributed
properties of the coefficients of $f(z)$ is a modification of
Ahlgren and Boylan's definition \cite{ab2} which is a special case
when $s=0$.

\begin{thm}
Suppose that $\lambda$ and $N$ are positive integers with $4\mid N$,
that $f(z)=\sum_{n=1}^\infty a(n)q^n\in
S_{\lambda+\frac{1}{2}}(\Gamma_0(N),\chi)\cap \mathbb{Z}[[q]]$ is a
half-integral weight cusp form, that $\chi$ is a real Dirichlet
character whose conductor divides $N$. Moreover, suppose that $M$ is
an odd positive integer. If for each $0\le r<M$, there is a positive
integer $n_r$ for which $a(n_r)\equiv r\pmod M,$ then the
coefficients of $f(z)$ are well-distributed modulo $M$.
\end{thm}

\begin{thm}
Suppose that $\lambda$ and $N$ are
positive integers with $4\mid N$, that $f(z)=\sum_{n=1}^\infty
a(n)q^n\in S_{\lambda+\frac{1}{2}}(\Gamma_0(N),\chi)\cap
\mathbb{Z}[[q]]$ is a half-integral weight cusp form, that $\chi$ is
a real Dirichlet character whose conductor divides $N$. Moreover,
suppose that $\ell$ is an odd prime and $j$ is a positive integer.
Then at least one of the following is true:

\noindent(1)The coefficients of $f(z)$ are well-distributed modulo
$\ell^j$.

\noindent(2)There are finitely many square-free integers $n_1, n_2,
\cdots, n_t$ for which
$$f(z)\equiv\sum_{i=1}^t\sum_{m=1}^\infty
a(n_im^2)q^{n_im^2}\pmod\ell.$$
\end{thm}

\noindent{\it Remark 1.2.} Under Ahlgren and Boylan's definition of
the coefficients of a half-integral weight modular form are well
distributed modulo $M$, Theorem 1.2 was proved by Bruinier and Ono
\cite{bo1,bo2} in the case when $j=1$ and by Ahlgren and Boylan
\cite{ab2} for any $j\ge1.$

As an application of Theorem 1.2, we consider the partition function
$p(n)$ and establish the following theorem.

\begin{thm}
If $\ell\ge5$ is prime and $j\ge1$ is an integer, then for any
integer $s\ge1$,
\begin{align*}
\sharp\{0\le n\le X\ |\ p(n)&\equiv r\pmod{\ell^j}\}\\
&\gg_{s,r,\ell^j}\left\{\begin{array}{ll}\frac{\sqrt{X}}{\log
X}(\log\log
X)^s \quad &\text{ if }r\not\equiv0\pmod{\ell^j},\\
X\quad&\text{ if }r\equiv0\pmod{\ell^j}.\end{array}\right.
\end{align*}
\end{thm}

There is an analogue result for integer weight modular forms.

\begin{thm} {\it Let $k\ge1$ and $M\ge1$ be integers with $M$ odd, $\chi$
be a real Dirichlet character,  and $f(z)=\sum_{n=1}^\infty
a(n)q^n\in S_{k}(\Gamma_0(N),\chi)\cap \mathbb{Z}[[q]]$ be an
integer weight cusp form with integral coefficients. If there exists
a coefficient $a(n_0)$ satisfying $(a(n_0),M)=1$, then for any
integer $ r$ and $s\ge 1$, }
$$\sharp\{1\le n\le X | a(n)\equiv r\pmod M\}
\gg_{f,s,r,M}\frac{{X}}{\log X}(\log\log X)^s.$$
\end{thm}
\noindent{\it Remark 1.3.} In the case $r\equiv0\pmod M$, Serre
\cite{s} proved that $\sharp\{n\le X\ |\ a(n)\equiv0\pmod M\}\gg_f
X$. (see also \cite[Theorem 2.65]{o2}).

Here we give an example as an application of Theorem 1.4.

\noindent{\bf Example.}  Let $a$ be a positive integer, $p$ an odd
prime. Recall that a $p^a$-regular partition of an integer $n$ is a
partition none of whose parts is divisible by $p^a$. We denote the
number of $p^a$-regular partitions of $n$ by $b_{p^a}(n)$. Penniston
\cite{p} proved that if $j$ is a positive integer and $1\le r\le
p^j$, then
$$\sharp\{n\le X\ |\ b_{p^a}(n)\equiv r(\text{ mod
}p^j)\}\gg\frac{X}{\log X}.$$By Proposition 2.2 of \cite{p}, for
every $j\ge3$ there exists an integer weight cusp form $F_j(z)$ such
that $F_j(z)\equiv \sum_{n=0}^\infty
b_{p^a}(n)q^{\frac{24n+p^a-1}{t}}(\text{ mod }p^j)$, where
$t=gcd(p^a-1,24)$. By Theorem 1.4, we can replace the bound above by
$\frac{{X}}{\log X}(\log\log X)^s$.

\section{Proofs of the results}

Let $p$ be a prime and $k$ be a positive integer. If
$f(z)=\sum_{n=1}^\infty a(n)q^n\in S_k(\Gamma_0(N),\chi),$ then the
action of the Hecke operator
$T|_{p,k,\chi}:S_k(\Gamma_0(N),\chi)\rightarrow
S_k(\Gamma_0(N),\chi) $ on $f(z)$ is defined by
\begin{align}f(z)|T_{p,k,\chi}:=\sum_{n=1}^\infty\left(a(pn)
+\chi(p)p^{k-1}a(n/p)\right)q^n.\end{align} If $p\nmid n$, then we
agree that $a(n/p)=0$.

Suppose that
$$F(z)=\sum_{n=1}^\infty a(n)q^n\in
S_{\lambda+\frac{1}{2}}(\Gamma_0(N),\chi).$$ Then the half-integral
weight Hecke operator $T_{p^2,\lambda,\chi}:
S_{\lambda+\frac{1}{2}}(\Gamma_0(N),\chi)\rightarrow
S_{\lambda+\frac{1}{2}}(\Gamma_0(N),\chi)$ on $F(z)$ is given by
\begin{align}F(z)|T_{p^2,\lambda,\chi}:=\sum_{n=1}^\infty\left(a(p^2n)
+\chi^\ast(p)\left(\frac{n}{p}\right)p^{\lambda-1}a(n)
+\chi^\ast(p^2)p^{2\lambda-1}a(n/p^2)\right)q^n,
\end{align}
where
$\chi^\ast$ is the Dirichlet character defined by
$\chi^\ast(n):=(\frac{(-1)^\lambda}{n})\chi(n)$, and $a(n/p^2)=0$ if
$p^2\nmid n$.

\begin{lem} Suppose that $M, k$ are positive integers,
and that $f(z)=\sum_{n=0}^\infty a(n)q^n\in S_k(\Gamma_0(N),\chi)$.
Then a positive proposition of the primes $p\equiv1\pmod{NM}$ have
the property that
$$f(z)|T_{p,k,\chi}\equiv2f(z)\pmod M.$$
\end{lem}

\begin{lem} {\it Suppose that $M, \lambda$ are positive integers and that
$f(z)=\sum_{n=0}^\infty a(n)q^n\in
S_{\lambda+\frac{1}{2}}(\Gamma_0(N),\chi)$, Then a positive
proposition of the primes $p\equiv1\pmod{NM}$ have the property
that}
$$f(z)|T_{p^2,\lambda,\chi}\equiv2f(z)\pmod M.$$
\end{lem}

Lemma 2.1 was observed by Serre \cite[Section 6.4]{s} for integer
weight cusp forms (or see \cite[Lemma 2.64]{o2}). By the Shimura's
theory on integer and half-integral weight modular forms, the result
holds for half-integral weight modular forms. For a proof of Lemma
2.2, see \cite[Lemma 2.2]{bo1}.

{\it Proof of Theorem} 1.1.  Since $f(z)\in
S_{\lambda+\frac{1}{2}}(\Gamma_0(N),\chi)\cap \mathbb{Z}[[q]]$, it
is clear that $f(z)\in
S_{\lambda+\frac{1}{2}}(\Gamma_0(2N\prod_rn_r),\chi)$. By Lemma 2.2,
the set of prime numbers
$$S(f,M):=\{p: p\equiv 1\pmod{2MN\prod_rn_r},
f(z)|T_{p^2,\lambda,\chi}\equiv2f(z)\pmod M\}$$ contains a positive
proportion of the primes. If $p\in S(f,M)$, then, by (3) for each
$n_r$
$$a(p^2n_r)+\chi^\ast(p)\left(\frac{n_r}{p}\right)p^{\lambda-1}a(n_r)
\equiv2a(n_r)\pmod M.$$ Since $p\equiv1\pmod {MN}$, we have
$$\chi^\ast(p)p^{\lambda-1}\equiv1\pmod M.$$
Assume that $n_r=2^{e(r)}\prod_jp_{j,r},$ where each $p_{j,r}$ is
odd. Then $p\in S(f,M)$ implies that $p\equiv1\pmod8$. This shows
that
$$\left(\frac{n_r}{p}\right)=\prod_j\left(\frac{p_{j,r}}{p}\right)
=\prod_j\left(\frac{p}{p_{j,r}}\right)=\prod_j\left(\frac{1}{p_{j,r}}\right)=1.$$
Therefore we conclude that for any $p\in S(f,M)$
$$a(p^2n_r)\equiv a(n_r)\pmod M\leqno{(4)}$$
and the coefficient of $q^{n_r}$ in $f(z)|T_{p^2,\lambda,\chi}\pmod
M$ is $a(p^2n_r)+a(n_r).$

For any $s$ distinct primes $p_1, p_2,\cdots, p_s\in S(f,M)$, we
claim that
$$a(p_s^2p_{s-1}^2\cdots p_1^2n_r)\equiv a(n_r)\pmod M.\leqno{(5)}$$ We
will prove the claim by induction. For any $p_1, p_2\in S(f,M)$,
$$(f(z)|T_{p_1^2,\lambda,\chi})|T_{p_2^2,\lambda,\chi}\equiv(2f(z))|T_{p_2^2,\lambda,\chi}
=4f(z)(\text{ mod }M).$$ Using (3) and comparing the coefficient of
$q^{n_r}$, we find that
$$a(p_2^2p_1^2n_r)+a(p_1^2n_r)+a(p_2^2n_r)+a(n_r)\equiv4a(n_r)(\text{ mod }M).$$
Combing (4) we obtain
$$a(p_2^2p_1^2n_r)\equiv a(n_r)(\text{ mod }M)$$ and the coefficient
of $q^{n_r}$ in
$(f(z)|T_{p_1^2,\lambda,\chi})|T_{p_2^2,\lambda,\chi}\pmod M$ is
$a(p_2^2p_1^2n_r)+3a(n_r).$ We assume that
$$a(p_{s-1}^2\cdots p_2^2p_1^2n_r)\equiv a(n_r)\pmod M\leqno{(6)}$$ and the coefficient
of $q^{n_r}$ in
$(f(z)|T_{p_1^2,\lambda,\chi})|T_{p_2^2,\lambda,\chi}\cdots
|T_{p_{s-1}^2,\lambda,\chi}(\text{ mod }M)$ is
$$a(p_{s-1}^2\cdots p_1^2n_r)+(2^{s-1}-1)a(n_r).$$
Then  Lemma 1.2 implies that
$$f(z)|T_{p_1^2,\lambda,\chi}|T_{p_2^2,\lambda,\chi}\cdots |T_{p_s^2,\lambda,\chi}
\equiv 2^sf(z)\pmod M.$$ Applying (3), we have
\begin{align*}a(p_s^2p_{s-1}^2\cdots p_1^2n_r)+(2^{s-1}-1)a(p_s^2n_r)
&+a(p_{s-1}^2\cdots p_1^2n_r)+(2^{s-1}-1)a(n_r)\\ &\equiv
2^sa(n_r)\pmod M.\end{align*} By (4) and (6), we prove (5).

It is clear that
$$\sharp\{1\le n\le X | a(n)\equiv r\pmod M\}\ge
\sharp\{p_s^2\cdots p_1^2n_r\le X\ |\ p_1,\cdots, p_s\in S(f,M)\}.$$
To estimate the right hand side, we use an argument as Landau. Let
$\delta$ denote the density of $S(f,M)$. We denote
$\pi_s\left(S(f,M),\cdots,S(f,M);X\right)$ by the number of
$n=p_1p_2\cdots p_s\le X$ with $p_i\in S(f,M), i=1,\cdots,s.$ Using
equation (2.28) of [15, Section 2.5], we have
\begin{align*}&\sharp\{p_s^2\cdots p_1^2n_r\le X\ |\ p_1,\cdots, p_s\in S(f,M)
\text{ are distinct }\}\\
&=\pi_s\left(S(f,M),\cdots,S(f,M);\sqrt{\frac{X}{n_r}}\right)\\
&\gg\frac{\delta^s}{(s-1)!}\frac{\sqrt{\frac{X}{n_r}}}{\log
\sqrt{\frac{X}{n_r}}}\left(\log\log
\sqrt{\frac{X}{n_r}}\right)^{s-1}\\
&\gg_{f,s,r,M}\frac{\sqrt{X}}{\log X}(\log\log X)^{s-1}.\end{align*}

Now we consider the case when $r\equiv0\pmod M$. By the Shimura's
theory and Serre's observation in \cite[Section 6.4]{s}, there is a
positive proportion of the prime $p\equiv-1\pmod {2MN\prod_rn_r}$
with the property that
$$f(z)|T_{p^2,\lambda,\chi}\equiv0\pmod M.$$ Fix one of such prime $p_0$.
Then it follows from (3) that for all $n$,
$$a(p^2_0n)\equiv-\chi^*(p)\left(\frac{n}{p_0}\right)p_0^{\lambda-1}a(n)-p_0^{2\lambda-1}
a(n/p_0^2)\pmod M.$$ Let $n=p_0m$ with $p_0\nmid m$. Then we find
that
$$a(p^2_0n)=a(p_0^3m)\equiv0\pmod M.$$
Since $m$ is arbitrary, we conclude that there is a positive
proportion of $n$ such that $a(n)\equiv0\pmod M$. \hfill$\Box$

{\it Proof of Theorem} 1.2. Suppose that $M$ is an odd positive
integer and $S(f,M)$ is as in the proof of Theorem 1.1. By Lemma 2.1
of \cite{bo1}, we find that if there exist $p_0\in S(f,M)$ and a
positive integer $n_0$ for which $(\frac{n_0}{p_0})=-1$, then for
every integer $r$, there exists an integer $n_r$ such that
$a(n_r)\equiv r\pmod M$. By Theorem 1.1, the coefficient of $f(z)$
are well-distributed modulo $M.$

Now let $M=\ell^j$. Suppose on the contrary that the coefficients of
$f(z)$ are not well-distributed modulo $\ell^j$. Then for every
$p\in S(f,\ell^j)$ and every integer $n$ with $a(n)\not\equiv0\pmod
{\ell^j}$, we have $(\frac{n}{p})\in\{0,1\}$. Following the argument
in the proof of Theorem 2.3 of \cite{bo1}, we conclude that there
are finitely many square-free integers $n_1, n_2,\cdots, n_t$ such
that
$$f(z)\equiv\sum_{i=1}^t\sum_{m=1}^\infty a(n_im^2)q^{n_im^2}\pmod\ell.\leqno{(7)}$$

\hfill$\Box$

{\it Proof of Theorem} 1.3. If $\ell=5,7,$ or $11$, then for every
integer $j\ge0$, by the construction outlined in \cite{ao}, we can
find positive integer $N_{\ell,j}, \lambda_{\ell,j}$, a quadratic
character $\chi_{\ell,j}$ and a modular form $G_{\ell,j}(z)\in
S_{\lambda_{\ell,j}+\frac{1}{2}}(\Gamma_0(N_{\ell,j}),\chi_{\ell,j})$
such that
$$G_{\ell,j}(z)\equiv\sum_{(\frac{-n}{p})=-1}p\left(\frac{
n+1}{24}\right)q^n\pmod {\ell^j}.$$ Suppose that $\ell\ge13$. Then
for every integer $j\ge0$,  Proposition 1 and 2 of \cite{a} imply
that there exists a modular form $F_{\ell,j}\in
S_{\frac{\ell^j-\ell^{j-1}-1}{2}}(\Gamma_0(576\ell),\chi_{12})$ for
which $$F_{\ell,j}\equiv\sum_{n=0}^\infty p\left(\frac{\ell
n+1}{24}\right)q^n\pmod {\ell^j}.$$ Following the argument in the
proof of Theorem 5 of \cite{ab2}, one can see that $G_{\ell,j}(z)$
and $F_{\ell,j}(z)$ can not have the form (7).  Theorem 1.3 follows
from Theorem 1.2 immediately. \hfill$\Box$

{\it Proof of Theorem} 1.4. Serre \cite[Section 6.4]{s} observed
that for any positive integer $M$, if $f(z)=\sum_{n=1}^\infty
a(n)q^n\in S_{k}(\Gamma_0(N),\chi)\cap \mathbb{Z}[[q]]$, then there
is a set of primes, say $V(M,f,i,n)$, of positive density with the
property that
$$a(np^i)\equiv(i+1)a(n)\pmod M$$for every pair of positive integers $i$ and $n$.
Therefore for any prime $\ell\in V(f,M,i_0,n_0)$, we have
$$a(n_0\ell^{i_0})\equiv(i_0+1)a(n_0)\pmod M.\leqno{(8)}$$
By Lemma 1.1, a positive proposition of the primes
$p\equiv1\pmod{NM}$ have the property that
$$f(z)|T_{p,k,\chi}=2f(z)\pmod M.$$
Denote the set of these primes by $T(f,M)$. For any $s$ distinct
primes $p_1, p_2, \cdots, p_s\in T(f,M)$, applying Lemma 1.1
repeatedly, we have
$$f(z)|T_{p_1,k,\chi}|T_{p_2,k,\chi}\cdots|T_{p_s,k,\chi}\equiv2^sf(z)\pmod M.$$
If $(p_1p_2\cdots p_s,\ell^{i_0}n_0)=1$, then by (2) we find that
$$a(p_sp_{s-1}\cdots p_1\ell^{i_0}n_0)\equiv 2^sa(\ell^{i_0}n_0)\pmod M.$$
Combing (8), we get
$$a(p_sp_{s-1}\cdots p_1\ell^{i_0}n_0)\equiv 2^s(i_0+1)a(n_0)\pmod M.$$
Since $(2a(n_0), M)=1$, as $i_0$ varies, $2^s(i_0+1)a(n_0)$ covers
all the residue classes modulo $M$.

For every $1\le r\le M-1$, choose $i_r$ such that
$2^s(i_r+1)a(n_0)\equiv r\pmod M$ and fix a prime $\ell_r\in
V(f,M,i_r,n_0)$. Let $\omega(\ell_r^{i_r}n_0)$ be the number of
distinct prime factors of $\ell_r^{i_r}n_0$ and $\delta_1$ be the
density of $T(f,M)$. Then by the argument of \cite[Section 2.5]{ni},
we have
\begin{align*}\sharp\{n\le X\ | a(n)&\equiv r\pmod M\}\\
&\ge\sharp\{p_s\cdots p_1\ell_r^{i_r}n_0\le X\ |\ p_1,\cdots, p_s\in
T(f,M)\text{ are distinct }\\ &\hskip5cm\text{ and }(p_1p_2\cdots p_s,\ell_r^{i_r}n_0)=1 \}\\
&\ge\pi_s\left(T(f,M),\cdots,T(f,M);{\frac{X}{\ell_r^{i_r}n_0}}\right)\\
&\hskip2cm-s\omega(\ell_r^{i_r}n_0)\pi_{s-1}\left(T(f,M),\cdots,T(f,M);{\frac{X}{\ell_r^{i_r}n_0}}\right)\\
&\gg\frac{\delta_1^s}{(s-1)!}\frac{{\frac{X}{\ell_r^{i_r}n_0}}}{\log
{\frac{X}{\ell_r^{i_r}n_0}}}\left(\log\log
{\frac{X}{\ell_r^{i_r}n_0}}\right)^{s-1}\\
&\gg_{f,s,r,M}\frac{{X}}{\log X}(\log\log X)^{s-1}.\end{align*} This
complete the proof of Theorem 4. \hfill$\Box$

\vskip.5cm

\begin{center}\textsc{Acknowledgements}\end{center}

The author thanks the National Natural Science Foundation of China
(Grant No.11026080) and the Natural Science Foundation of Education
Department of Henan Province (Grant No. 2009A110001).

\end{document}